\documentclass[11pt,dvips]{amsart}
\usepackage{amscd,amsthm,amsmath,amsgen,amsfonts, amssymb, latexsym, epsfig, euler, epic, graphicx, psfrag}

\numberwithin{equation}{section}
\numberwithin{figure}{section}
\usepackage{amscd}

\input xy
\xyoption{all}

\newtheorem{theorem}{Theorem}[section]
\newtheorem{corollary}[theorem]{Corollary}
\newtheorem{lemma}[theorem]{Lemma}
\newtheorem{prop}[theorem]{Proposition}
\newtheorem{definition}[theorem]{Definition}
\newtheorem{remark}[theorem]{Remark}

\newtheorem{example}[theorem]{Example}

\def\Z{\mathbb{Z}}

\def\Q{\mathbb{Q}}
\def\R{\mathbb{R}}
\def\C{\mathbb{C}}
\def\CP{\mathbb{CP}}

\def\t{\mathfrak{t}}

\def\L{\mathcal{L}}

\def\onestep{E}
\def\GKMedges{E_{GKM}}
\def\spaceafterx{\!\!\!\!\!\!\!}
\def\spacebeforeeta{\!\!\!\!\!\!\!\!}

\newcommand{\Stab}{\mathrm{Stab}}
\newcommand{\excise}[1]{}

\begin{document}
\title{Towards generalizing Schubert calculus in the symplectic category}
\author{R. F. Goldin}\address{rgoldin@math.gmu.edu\\ George Mason University\\ MS 3F2, 4400 University Dr.\\ Fairfax, VA 22030}
\author{S. Tolman}\address{stolman@math.uiuc.edu\\ University of Illinois Urbana-Champaign\\ 1409 W. Green Street\\ Urbana, IL 61801}
\thanks{The first author  was partially supported by NSF-DMS Grant \#0606869,
and the second author by  NSF-DMS Grant \#0707122.}
\begin{abstract}
The main purpose of
this article is to extend some of the ideas from Schubert calculus to the more
general setting of Hamiltonian torus actions on compact symplectic manifolds
with isolated fixed points.
Given a generic component $\Psi$ of the moment map, which is a Morse function,
we define a canonical class $\alpha_p$ in
the equivariant cohomology of the manifold $M$ for each fixed point $p \in M$.
When they exist, canonical classes form a natural basis of the equivariant
cohomology of $M$; in particular, when $M$ is a flag variety, these classes are
the equivariant Schubert classes. We show that the restriction of
a canonical class $\alpha_p$ to a fixed point $q$
can be calculated by a rational function
which depends only on the value of the moment map,
and the restriction of other canonical classes to points of index exactly two higher.
Therefore, the structure constants can be calculated by a similar rational function.
Our restriction formula is {\em manifestly positive} in many cases, including when $M$
is a flag manifold. Finally, we prove the existence of integral
canonical classes in the case that $M$ is a GKM manifold and $\Psi$ is {\em index
increasing}. In this case, our restriction formula specializes to an easily
computable rational sum which depends only on the GKM graph.
\end{abstract}
\maketitle

\tableofcontents

\section{Introduction}\label{se:introduction}

Let $T$ be a compact torus with Lie algebra $\t$ and
lattice $\ell \subset \t$.
Suppose that 
$T$ acts on a compact symplectic manifold $(M,\omega)$
with isolated fixed points
and moment map $\Phi \colon M \to \t^*$,
where $\t^*$ is  dual  to $\t$.
Then
$$
\iota_{X_\xi}\omega = -d \Phi^\xi \quad \forall \  \xi  \in \t,
$$
where  
$X_\xi$ denotes  the vector field on $M$ generated by the action 
and
$\Phi^\xi \colon M \to \R$ is defined by 
$\Phi^\xi(x) = \langle \Phi(x), \xi \rangle$. 
Here, $\langle \cdot, \cdot \rangle$ is the natural pairing
between $\t^*$ and $\t$.

If $\xi \in \t$ is {\bf generic}, that is,
if $\langle \eta, \xi \rangle \neq 0$ for each weight $\eta \in \ell^* \subset \t^*$
in the symplectic  representation $T_p M$ for  every $p$ in the fixed set $M^T$,
then $\Psi = \Phi^\xi \colon M \to \R$ is a Morse function with critical set $M^T$.
Given $p \in M^T$, the
negative tangent bundle $\nu^-(p)$ is a  representation
with no fixed sub-bundle.
Hence, the index of $\Psi$ at $p$ is even;
let $\lambda(p)$ denote  half the index of $\Psi$ at $p$.
The individual weights of this
representation are 
well defined and non-zero;
our convention for the moment map implies that these weights
are exactly the {\bf positive  weights}  of the $T$ action on $T_p M$,
that is, the weights $\eta$ such that $\langle \eta, \xi \rangle > 0$.
Let $\Lambda_p^-$ denote the product of these weights.
(Conversely, the weights in the positive tangent bundle are the
{\em  negative weights} of the $T$ action on $T_pM$.)  
Finally, for all $q \in M^T$ the
inclusion $q\hookrightarrow M$  induces a map 
$H_T^*(M)\rightarrow H_T^*(q)$ in equivariant cohomology;
let $\alpha(q)$ denote
the image of a class $\alpha \in H_T^*(M)$ under this map.

\begin{definition}\label{de:canonical}
Let a torus $T$ act on a compact symplectic manifold $(M,\omega)$
with isolated fixed points and moment map $\Phi \colon M \to \mathfrak{t}^*$.
Let $\Psi = \Phi^\xi :M\rightarrow\R$ 
be a generic component of the moment map.
A cohomology class $\alpha_p \in H^{2 \lambda(p)}_T(M;\Q)$ is
the {\bf canonical class} at a fixed point $p$ with respect to $\Psi$ if
\begin{enumerate}
\item[$(1)$] $\alpha_p(p) = \Lambda_p^- $
\item[$(2)$] $\alpha_p(q) = 0$ for all $q \in M^T \smallsetminus \{p\}$ such that
$\lambda(q) \leq \lambda(p)$.\footnote{Note that $(2)$ is stronger than the frequently encountered condition that $\alpha_p(q)=0$ for all $q\in M^T\smallsetminus \{p\}$ such that $\Psi(q)\leq \Psi(p)$. See Lemmas~\ref{le:pclass} and~\ref{le:2prime}.}
\end{enumerate}
Moreover, we say that the canonical class
$\alpha_p$ is {\bf integral} if  
$\alpha_p \in H^{2 \lambda(p)}_T(M;\Z).$\footnote{
Since the fixed points are isolated, $H^*_T(M;\Z)$ is torsion free;
see Lemma~\ref{le:pclass}.  Therefore, we can naturally
identify $H^*_T(M;\Z)$ with a subgroup of $H^*_T(M;\Q)$.
} 
\end{definition}

We cannot always find canonical classes; see Example~\ref{ex:CP2}. 
However,  
each  canonical class is unique and can be thought of as an 
equivariant Poincar\'e dual to the  
closure of the stable manifold. 
If $\alpha_p$ exists for all $p \in M^T$, then $\{\alpha_p\}$ forms a basis of $H_T^*(M)$ as a module over $H^*(BT)$.
Since the fixed set is isolated, the natural 
restriction map $H_T^*(M;\Z) \to H^*(M;\Z)$ is surjective; 
under this map,
the canonical classes also define a basis for the ordinary cohomology  $H^*(M)$.
In the case that $M = G/B$, 
where  $G$ is a complex semi-simple Lie group (of any type) 
and $B$ is a Borel subgroup,
the equivariant Schubert classes are  canonical classes. 
Under the map to ordinary cohomology, they are 
exactly the Poincar\'e duals to Schubert varieties in ordinary cohomology.
Hence, our work is a direct generalization of that setting.

This paper is concerned with a new formula
for how to restrict canonical cohomology classes
to  fixed points. 
Since the fixed points are isolated,
the inclusion of the fixed point set $M^T$ into $M$ induces an injection 
$H_T^*(M;\Z)\rightarrow H_T^*(M^T;\Z)$, 
where the latter ring is a direct sum of polynomials rings. 
Thus each cohomology class on $M$ may be described by an integral polynomial associated to
each fixed point. 
Once the restriction of canonical classes is known at each fixed point, one 
can easily derive a formula for the structure constants in the (equivariant) 
cohomology ring. (See \cite{GZ:Thomclasses}.) 
Recall that the structure constants for $H_T^*(M)$ are the set $c_{pq}^r\in H_T^*(M)$ given by
$$
\alpha_p\alpha_q = \sum_{r\in M^T} c_{pq}^r \alpha_r.
$$
Conversely, the structure constants also provide a formula for the restrictions.

Our formulas have  some 
echoes in the literature; S. Billey \cite{Bi:Kostantpolynomials} found a 
different manifestly positive formula for the restriction of equivariant 
Schubert classes when $M=G/B$. 
V. Guillemin and C. Zara \cite{GZ:Thomclasses} found a non-positive path formula for the restrictions in the case of GKM graphs, which we discuss in more detail below.

Our main contribution in this article can be seen as an inductive formula
for the restriction of canonical classes to fixed points; we prove this in 
Section~\ref{se:induction}. The formula depends on 
only
the values of 
the moment map and 
$\alpha_r(r')$, where $r$ and $r'$ are fixed points whose indices differ by two.

Given a directed graph with vertex set $V$ and edge set $E \subset
V \times V$, a {\bf path}
from a vertex $p$ to a vertex $q$ is a $(k+1)$-tuple
${\bf r} = (r_0,\ldots,r_k) \in V^{k+1}$ so that $r_0 = p$, $r_k = q$,
and $(r_{i-1},r_i) \in E$ for all $1 \leq i \leq k$;
let $|{\bf r} | = k$ denote the {\bf length} of ${\bf r}$.

\begin{theorem}\label{th:pathformula}
Let a torus $T$ act on a compact symplectic manifold
$(M,\omega)$ with isolated fixed points and moment map $\Phi \colon M \to \t^*$.
Let $\Psi = \Phi^\xi$ be a generic component of the moment map.
Assume that there exists a canonical class
$\alpha_p \in H_T^{2 \lambda(p)} (M;\Q)$ for  all $p \in M^T$.

Define an oriented graph with vertex set $V = M^T$ and edge set
$$
\onestep = \left\{ \left. (r,r') \in M^T \times M^T \right| \lambda(r') - \lambda(r) = 1 \mbox{ and } \alpha_r(r') \neq 0
\right\}.
$$
Given $p$ and $q$  in $M^T$, let $\Sigma_p^q$ denote the set of paths from $p$ to $q$ in $(V,E)$; then
\begin{equation}\label{eq:generalrestriction}
\alpha_p(q) =
 \Lambda_q^-
\sum_{{\bf r} \in \Sigma_p^q} \
\prod_{i=1}^{|{\bf r}|}
\frac {\Phi(r_{i}) - \Phi(r_{i-1})} {\Phi(q)-\Phi(r_{i-1})}
\cdot 
\frac{\alpha_{r_{i-1}}(r_i)}{ \Lambda^-_{r_i}}.
\end{equation}
\end{theorem}

\begin{remark} [{\bf Positivity}] \label{positivity} \rm
We say that $\alpha \in S(\t^*)$ is {\bf positive} if $\alpha(\xi) > 0$
and {\bf negative} if $\alpha(\xi) < 0$.
In some cases, the restriction $\alpha_p(q)$ is itself negative; 
see Example~\ref{ex:nonkahler}.
A fortiori, in these examples some of the summands in 
\eqref{eq:generalrestriction}
are negative.
However, whenever $\alpha_p(q) \geq 0$ for all $p$ and $q \in M^T$
such that $\lambda(q) = \lambda(p) + 1$,
our formula is {\bf manifestly positive}, in the sense that
each summand is positive.
To see this, note that $\Lambda_q^-$ 
and $\Lambda_{r_i}^-$
are positive by definition,
$\Phi(r_i) - \Phi(r_{i-1})$ and $\Phi(q) - \Phi(r_{i-1})$ are
positive by Corollary~\ref{co:increasing},
and 
$ \alpha_{r_{i-1}}(r_i) $ 
is positive by assumption. 

For example, for flag varieties $G/B$ of semisimple Lie groups
the canonical classes are Schubert classes; see \cite{BGG}. 
In this case, the restriction $\alpha_p(q)$ is positive for
all $p$ and $q$ 
by \cite{Bi:Kostantpolynomials}.
Alternatively, it is very easy to check this directly
when  $\lambda(q) = \lambda(p) + 1$;
see Section~\ref{se:examples} for the case
$G=Sl(n,\C)$  and \cite{ST} for the general case.
\end{remark}

\begin{corollary}
Consider the situation described in Theorem~\ref{th:pathformula}.  If there is no path in 
$(V,E)$ from a fixed point $p$ to a fixed point $q$, then $\alpha_p(q) =  0$.
Moreover, if
$\alpha_p(q) \geq 0$ for all  $p$ and $q$ in $M^T$ such that $\lambda(q)
= \lambda(p) + 1$, then $\alpha_p(q) \geq 0$ for all $p$ and $q$ in $M^T$ and
$\alpha_p(q) > 0$ exactly if there is at least one
path from $p$ to $q$.
\end{corollary}

We now restrict our attention to an important special case
where it is especially easy to make these calculations:
GKM spaces.
Let a torus $T$ act on a compact symplectic manifold $(M,\omega)$ with
moment map $\Phi:M\rightarrow \mathfrak{t}^*$.
We say that $(M,\omega,\Phi)$ is a {\bf GKM space}
if $M$ has isolated fixed points and if,
for every codimension one subgroup $K\subset T$,  
every connected component of the fixed submanifold
$M^K$ has dimension two or less.
\begin{definition}\label{GKMgraph}
Let $(M,\omega,\Phi)$ be a GKM space. We define the {\bf GKM graph} to be the
labelled directed graph $(V,\GKMedges)$ given as follows.
The vertex set $V$ is the fixed set $M^T$;
we label each $p \in M^T$ by its moment image $\Phi(p) \in \t^*$.
The edge set $\GKMedges$ consists of
pairs of distinct points $(p,q) \in V \times V$ such that
there exists a codimension one subgroup $K\subset T$ so that
$p$ and $q$ are contained in the same component $N$ of $M^K$.
We label the edge $(p,q)$ by the weight $\eta(p,q)\in \ell^*$
associated to the representation of $T$ on $T_qN\simeq \C$.
\end{definition}

Let $\Psi= \Phi^\xi$ be a generic component of the moment map. 
Note that $\lambda(p)$ is the number of edges  $(r,p) \in \GKMedges$ such that 
 $\Psi(r) < \Psi(p)$; moreover, $\Lambda_p^- = \prod \eta(r,p)$,
where the product is over all such edges.
We say that $\Psi$ is {\bf index increasing} 
if $\Psi(p)<\Psi(q)$ implies 
that $\lambda(p)<\lambda(q)$
for all $(p,q) \in \GKMedges$. See Example~\ref{ex:indexincreasing} and 
Remark~\ref{rm:increasing}.

Given any weight $\eta \in \ell^*$, 
the projection which takes
$X \in \t^*$ to $ X-\frac{\langle X,\xi\rangle}{\langle\eta,\xi\rangle}\eta 
\in \xi^\perp \subset \t^*$ 
naturally induces a endomorphism $\rho_\eta$ of  $S(\t^*)$, the
symmetric algebra on $\t^*$.
Since $M$ is a GKM space,
the weights at each fixed point  are pairwise linearly independent; hence,
$\rho_{\eta(r,r')}(\Lambda^-_r) \neq 0$ and 
$\rho_{\eta(r,r')} \left(\frac{\Lambda^-_{r'}}{\eta(r,r')}\right) \neq 0$
for all $(r,r') \in \GKMedges$.
Following \cite{GZ:Thomclasses}, 
we define
\begin{equation}\label{eq:Theta}
\Theta(r,r') = \frac{\rho_{\eta(r,r')} \left(\Lambda^-_r \right)}
{\rho_{\eta(r,r')} \left(\frac{\Lambda^-_{r'}}{\eta(r,r')}\right)} 
\in S(\t^*)_0 \quad \forall \ (r,r') \in \GKMedges,
\end{equation}
where $S(\t^*)_0$ is the field of fractions of $S(\t^*)$.

\begin{theorem}\label{th:GKMpathformula}
Let $(M,\omega,\Phi)$ be  GKM space. Let 
$\Psi = \Phi^\xi$ be
a generic component
of the moment map; assume 
that $\Psi$ is index increasing.
Define an oriented graph with vertex set $V= M^T$
and edge set
$$
\onestep= \left\{(r,r')\in \GKMedges\ \left|
\ \lambda(r')-\lambda(r)=1 \right. \right\},
$$
where $(V,\GKMedges)$ is the GKM graph associated to $M$.
Then
\begin{itemize}
\item There exists a canonical class
$\alpha_p\in H_T^{\lambda(p)}(M;\Z)$ for all $p \in M^T$.
\item
Given $p$ and $q$ in $M^T$, 
let $\Sigma_p^q$ denote the set of paths from $p$ to $q$ in $(V,\onestep)$; then
\begin{equation}\label{GKMsum}
\alpha_p(q) =
 \Lambda_q^-
\sum_{{\bf r} \in \Sigma_p^q} \
 \prod_{i=1}^{|r|}
\frac {\Phi(r_{i}) - \Phi(r_{i-1})} {\Phi(q)-\Phi(r_{i-1})}\cdot 
\frac{\Theta(r_{i-1},r_i)}{\eta(r_{i-1},r_{i})}. 
\end{equation}
\item $\Theta(r,r')\in \Z \smallsetminus \{0\}$ for all $(r,r')\in E$. \end{itemize}
\end{theorem}

\begin{remark}\rm \label{re:equivalenceTheta}
A straightforward calculation shows that, since $\Theta(r,r')$ is an integer,
\begin{equation}\label{eq:Thetasatisfies}
\frac{\Lambda_{r'}^-}{\eta(r,r')}\Theta(r,r') = \Lambda_r^- \mod \eta(r,r').
\end{equation}
Moreover, since 
$\frac{\Lambda_{r'}^-}{\eta(r,r')}$ is not a multiple of $\eta(r,r')$,
equation~\eqref{eq:Thetasatisfies} has a unique solution
and so provides an alternative definition of $\Theta$.
\end{remark}

\begin{remark} \rm
In fact,
since $\Phi(r') - \Phi(r)$ is a positive multiple of
$\eta(r,r')$ for all $(r,r') \in \GKMedges$,
the formula \eqref{GKMsum} is a manifestly positive   
exactly if $\Theta(r,r') >0$  for all $(r,r') \in \GKMedges$;
cf. Remark~\ref{positivity}.
However, $\Theta(r,r')$ is not always positive; 
see Example~\ref{ex:nonkahler}.
\end{remark}

\begin{example}\rm \label{cpn}
$e_1,\ldots,e_{n+1}$ denote  the standard basis for $\R^{n+1}$, and 
$x_1,\ldots,x_{n+1}$ denote the dual basis for $(\R^{n+1})^*$.
Let  $S^1_\Delta \subset (S^1)^{n+1}$ be the diagonal circle, and
$T = (S^1)^{n+1}/S^1_\Delta$
the quotient torus.
The standard action of $(S^1)^{n+1}$ on $\C^{n+1}$ induces a symplectic  action of $T$
on complex projective space
$(\CP^n,\omega)$
with  moment map  $\Phi \colon \C P^n \to \t^*$, where
$\t^* = \left\{ \left. \sum 
a_i x_i \in  (\R^{n+1})^* \ \right| \ \sum a_i = 0 \right\},$
and
 $$\Phi([z_1:\dots: z_{n+1}]) =  \sum_{i=1}^{n+1} \left( \frac{1}{n+1} - \frac{ |z_i|^2 }{ \sum_j |z_j|^2 }  \right)   x_i.$$
It is straightforward to check that $(\C P^n, \omega, \Phi)$ is a GKM space,
and that the associated GKM graph
is the complete directed graph on  $n+1$ vertices 
$p_1,\ldots,p_{n+1}$, where 
$p_i = [e_i \otimes \C]$.
Moreover,
$\Phi(p_i) = \frac{1}{n+1}\sum_j x_j - x_i$ for all $i$,
and $\eta(p_i,p_j) = x_i - x_j$ for all $i \neq j$.

Let $\Psi = \Phi^\xi$, where 
$\xi = -(0, 1,2,\dots,n)$.
Then $\Psi(p_i) - \Psi(p_j) = \langle - x_i + x_j, \xi \rangle = i - j$ 
is positive 
exactly if  $j < i$,
and so
$\Lambda_{p_i}^- = \prod_{j < i}( x_j -x_i)$ and 
$\lambda(p_i) = i - 1$ for all $i$.
Therefore, $\Psi$ is index increasing 
and $E = \{ (p_1,p_2), \ldots, (p_{n},p_{n+1}) \}$.
In particular, there is exactly one path $(p_i,p_{i+1},\ldots,p_j)$ from
$p_i$ to $p_j$ in $E$ if $i \leq j$; otherwise, there is none.
Finally,  since $\rho_{x_{l-1}-x_l}(x_l) = \rho_{x_{l-1}-x_l}(x_{l-1})$,
\begin{equation*}
\Theta(p_{l-1},p_l) = \frac
{\rho_{\eta(p_{l-1},p_l)} \left(\Lambda^-_{p_{l-1}} \right)}
{\rho_{\eta(p_{l-1},p_l)} \left(\frac{\Lambda^-_{p_l}}{\eta(p_{l-1},p_l)}\right)} 
=\frac
{\rho_{x_{l-1} - x_l} \left( \prod_{m<l-1}  (x_m-x_{l-1}) \right) }
{\rho_{x_{l-1} - x_l}\left(\prod_{m<l-1}(x_{m}-x_l) \right)} 
= 1 \quad \forall \ l.
\end{equation*}
Thus, by Theorem~\ref{th:GKMpathformula},
\begin{align*}
\alpha_{p_i}(p_j) &= 
\Lambda_{p_j}^- \prod_{l=i+1}^j 
\frac{\Phi(p_l)-\Phi(p_{l-1})}{\Phi(p_j)-\Phi(p_{l-1})}\cdot 
\frac{\Theta(p_{l-1},p_l)}{\eta(p_{l-1},p_{l})}  \\
&= \prod_{l=1}^{j-1} (x_l-x_j) \prod_{l=i+1}^j \frac{x_{l-1}-x_{l}}{x_{l-1}-x_{j}} \cdot 
\frac{1}{x_{l-1}-x_{l}} \\
&=\prod_{l=1}^{j-1}(x_l-x_j) \prod_{l=i}^{j-1}\frac{1}{x_{l}-x_j}
=\prod_{l=1}^{i-1} (x_l-x_j).
\end{align*}
\end{example}

\begin{remark}\rm
Guillemin and Zara also give a formula for $\alpha_p(q)$
for GKM spaces as a sum over paths in \cite{GZ:Thomclasses}.
In fact, their formula is identical to ours in
the case that $\lambda(q) - \lambda(p) = 1$, and
also works in a slightly broader context.
However, in general the formulas are quite different.
For example, their formula for $\alpha_p(q)$ includes a contribution for each
path $r$ from $p$ to $q$ in $(V,\GKMedges)$ such that
$\Psi(r_i) > \Psi(r_{i-1})$ for each edge $(r_{i-1},r_i)$;
our formula only includes a contribution from
a subset of such paths -- those such that $\lambda(r_i) = \lambda(r_{i-1}) + 1$.
In practice, this means that we sum over many fewer paths. For example,
if $M = \C P^n$ their formula for $\alpha_p(q)$ contains $2^{\lambda(q) - \lambda(p) - 1}$ terms,
whereas ours contains just one term; see Example~\ref{cpn}.
Moreover, their formula is almost never  manifestly positive, in the
sense  described above.
\end{remark}

In \cite{Kn:DHmeasure}, Knutson gives a positive formula for the
Duistermaat-Heckman measure of a torus action
on a smooth  algebraic variety\footnote{
The restriction to the algebraic category is only implicit.}  with an
invariant Palais-Smale metric, and suggests a technique for computing
the Duistermaat-Heckman  measure of certain subvarieties.
In fact, in the case that $M$ is an algebraic variety and
there exists an invariant Palais-Smale metric, it is
possible to use the results of \cite{Kn:DHmeasure} to give an alternate
proof of Theorem~\ref{th:pathformula}.  We hope to do this in our next paper;
we also plan to use Theorem~\ref{th:pathformula}  to extend his formula
for the Duistermaat-Heckman measure to the non-algebraic case.
However, in this greater generality, the summands in the formula
are {\em not} always positive.  This occurs, for example,
in the manifold considered in  Example~\ref{ex:nonkahler}.
After we initially announced these results, Knutson showed
that he could extend his results in \cite{Kn:DHmeasure} by
dropping the condition that there exist an invariant Palais-Smale
metric; see \cite{Kn:Localization}.

Finally, several techniques have recently been discovered which use the ideas
in this paper
to find a positive integral formula in certain important cases, including when $M$ is a flag manifold (\cite{ST},\cite{Za:positivity}).

We would like to thank Victor Guillemin, whose questions inspired this project.
We would also like to  thank Sara Billey, Allen Knutson, Catalin Zara, and Silvia Sabatini
for many helpful discussions.

\section{Canonical Classes}\label{se:canonical}

In this section, we demonstrate
some properties of canonical classes.
In particular,  we show that if  they  exist then they form a natural basis for
$H_T^*(M;\Z)$ as a $H^*(BT;\Z)$ module.
Additionally,   they do exist in a number of important cases.

For this purpose, it  is natural to work in a slightly more
general context.  Let a torus $T$ act on a compact
oriented manifold $M$ with isolated fixed points.
An invariant Morse function $\Psi \colon M\rightarrow \R$ is a
{\bf formal moment map}\footnote{In \cite{GGK}, these are
called ``non-degenerate abstract moment maps.''}
if the critical set of $\Psi$  is exactly the fixed point set
$M^T$.
As we saw in the introduction, if $M$ is symplectic and
the action  is  Hamiltonian, then any generic component of the moment map is a formal moment map.
The cohomological properties of symplectic manifolds
with Hamiltonian actions described in the introduction
continue to hold in the more
general case of formal moment maps; see Appendix G of \cite{GGK}.
In particular,
the restriction map from
$H^*_T(M;\Z)$ to $H^*_T(M^T;\Z)$  is injective in
this case.

Let $p$ be a critical point for a formal moment map $\Psi: M \to \R$.
Since $\nu^-(p)$ is a real representation with
no fixed subbundle, the index of $\Psi$ at $p$ is even;
let $\lambda(p)$ denote  half the index of $\Psi$ at $p$.
The signs of the individual weights of this
representation are not well defined. However,
if we fix an orientation on the negative normal bundle $\nu^-(p)$
then  the product of the weights  -- which we will denote $\Lambda_p^-$ --
is well-defined.

\begin{definition}\label{de:canclass}
Let a torus $T$ act on a compact oriented manifold $M$ with isolated fixed points,
 and let $\Psi \colon M \to \R$ be a formal moment map.
We say that a  cohomology class $\alpha_p \in H^{2 \lambda(p)}_T(M;\Q)$ is
a {\bf canonical class}  (with respect to $\Psi$) at a fixed point $p$ if there exists an orientation on $\nu^-(p)$ such that
\begin{enumerate}
\item[$(1)$] $\alpha_p(p) = \Lambda_p^- $, and
\item[$(2)$] $\alpha_p(q) = 0$ for all $q \in M^T \smallsetminus \{p\}$ such that
$\lambda(q) \leq \lambda(p)$.
\end{enumerate}
Moreover,
the canonical class is {\bf integral} if  $\alpha_p \in H^{2 \lambda(p)}_T(M;\Z)$.
\end{definition}

Canonical classes do not always exist.

\begin{example}\label{ex:CP2}\rm
Let the 
torus $ T= (S^1)^2$ act on $\CP^2$ by
$$ (t_1, t_2) \cdot [z_1:z_2:z_3]
= [t_1 z_1: t_2 z_2: z_3].$$
Let $M$
be the blow-up of  $\CP^2$
at $[0:0:1]$, and let 
$\Phi \colon M \to \R^2$ be the moment map for the
induced $T$ action.
Let $\Psi = \Phi^\xi,$
where $\xi = (1,-1) \in \R^2$.
Label the four fixed points
$p_1,\ldots,p_4$ so that
$\Psi(p_1) < \dots < \Psi(p_4)$.
There exists
a basis  $\gamma_1,\dots,\gamma_4$ for $H_T^*(M;\Z)$  as a $H^*(BT;\Z)$ module so that
$\gamma_i(p_j) = 0 \ \forall \ j < i$ and
\begin{gather*}
\gamma_1(p_1) =
\gamma_1(p_2) =
\gamma_1(p_3) =
\gamma_1(p_4) = 1;  \quad
\gamma_2(p_2) = \gamma_2(p_3) = x_1, \ \  \gamma_2(p_4) = x_1 - x_2; \\
\gamma_3(p_3) = \gamma_3(p_4) = x_1 - x_2; \quad
 \mbox{and}\quad
\gamma_4(p_4) = -x_2 (x_1 - x_2).
\end{gather*}
A straightforward calculation shows that there  is no canonical class for $p_2$.
(Although this example is a GKM space,
it is consistent with  Theorem~\ref{th:GKMpathformula}
because $\Psi$ is not index increasing.)
\end{example}

However, if canonical classes do exist, they give a natural basis for $H_T^*(M)$.

\begin{prop}\label{prop:basis}
Let a torus $T$ act on a compact oriented manifold $M$ with isolated fixed
points,  and let $\Psi \colon  M \to \R$ be a formal moment map.
Fix an orientation on $\nu^-(p)$ for all $p \in M^T$.
If there exists a canonical class $\alpha_p \in H_T^*(M;\Q)$
for all $p \in M^T$, then
the  classes $\{\alpha_p\}_{p \in M^T}$
are a natural  basis for  $H_T^*(M;\Q)$
as a module over $H^*(BT;\Q)$. 
Moreover, if the canonical classes are integral then
they
are a natural  basis for  $H_T^*(M;\Z)$
as a module over $H^*(BT;\Z)$.
\end{prop}

\begin{proof}
This is an immediate consequence of  Lemmas~\ref{le:pclass},
\ref{le:canprop}, and \ref{le:2prime} below.
\end{proof}

\begin{remark}\rm\label{sensitive}
This basis does not depend very sensitively on the choices that we have
made; it only depends on the $\Lambda_p^-$'s at each fixed point.
For example, let a torus $T$ act on a compact symplectic manifold $(M,\omega)$
with isolated fixed points and moment map $\Phi \colon M \to \mathfrak{t}^*$.
Let $\Psi = \Phi^\xi$  be a generic component of the moment map.
Let $\alpha_p \in H^{2\lambda(p)}_T(M;\Z)$ be a canonical class at $p$ for $\Psi$.
If $\Psi' = \Phi^{\xi'}$ is another generic component,
then $\alpha_p$ is also
a canonical class for $\Psi'$, as long as
$$ \langle \xi, \eta \rangle > 0 \leftrightarrow
\langle \xi', \eta \rangle > 0   \quad \forall \ \eta \in \Pi_p \mbox{ and }
p \in M^T,$$
where $\Pi_p$ denotes  the set of weights at $p$.
Similarly, if  $\omega' \in \Omega^2(M)$ is another invariant symplectic
form
with moment map $\Phi' \colon M \to \t$,
then $\alpha_p$ is also a canonical class for $\Psi' = (\Phi')^\xi$ as long
as $\omega$ and $\omega'$ are deformation equivalent.
\end{remark}

The following lemma is a key ingredient in our proof
that canonical classes exist in certain cases.
In particular, the existence of these closely related classes is guaranteed by
straightforward Morse theoretic arguments.

\begin{lemma} [Kirwan] \label{le:pclass}
Let a torus $T$ act on a compact oriented manifold
$M$ with isolated fixed points,
and let $\Psi: M \to \R$ be a formal moment map.
For every fixed point $p$
and  for each orientation on $\nu^-(p)$,
there exists
an integral cohomology  class $\gamma_p \in H_{T}^{2\lambda(p)}(M;\Z)$
so that
\begin{enumerate}
\item[$(1)$]  $\gamma_p(p) =\Lambda^-_p $, and
\item[$(2^\prime)$] $\gamma_p(q) = 0$ for every
$q \in M^T \smallsetminus \{p\}$ such that
$\Psi(q)\leq \Psi(p)$.
\end{enumerate}
Moreover, for any such classes, the
$\{ \gamma_p \}_{p \in M^T}$ are a basis
for  $H_{T}^*(M;\Z)$ as a module over
$H^*(BT;\Z)$.
\end{lemma}

Kirwan proved this result for
rational cohomology classes on compact symplectic manifolds with
Hamiltonian actions \cite{Ki:quotients}
(see also \cite{TW:injectivity}). The proof generalizes easily to the
case of formal moment maps, and to
integral cohomology when the fixed points are isolated.

\begin{corollary}\label{co:towardindexstatement}
Let a torus $T$ act on a compact oriented manifold
$M$ with isolated fixed points,  and let $\Psi: M \to \R$
be a formal moment map.
Given a point  $p \in M^T$ and a class $\beta \in H^{2 i}_T(M;\Q)$
such that
 $\beta(q) = 0$  whenever $q \in M^T$ satisfies $\Psi(q) < \Psi(p)$, 
\begin{itemize}
\item The restriction $\beta(p)= x\Lambda_p^-,$ where $x\in H^*(BT;\Q)$.
\item In particular, if $\lambda(p) > i$
then $\beta(p) = 0$.
\item If $\beta\in H_T^*(M;\Z)$ is integral 
then $\beta(p) = x \Lambda_p^-$, where $x\in H^*(BT;\Z)$ is integral.
\end{itemize}
\end{corollary}
\begin{proof}
By Lemma~\ref{le:pclass}
for each $q \in M^T$ we can fix an orientation on
$\nu^-(q)$ and choose a class $\gamma_q \in H^*_T(M;\Z)$
which satisfies properties $(1)$ and $(2^\prime)$.
Moreover, we can write
$\beta = \sum_{q \in M^T}  x_q \gamma_q$,
where $x_q$ lies in $H^*(BT;\Q)$
for all $q$.
If $\beta$ is integral, then each $x_q$ lies in  $H^*(BT;\Z)$.

If $x_q = 0$ for all $q \in M^T$ so that $\Psi(q) < \Psi(p)$,
then properties $(1)$ and $(2^\prime)$ together imply
that $\beta(p) = x_p \Lambda^-_p$.
Otherwise, there  exists $q \in M^T$ so that  $\Psi(q) < \Psi(p)$
and $x_q \neq 0$, but
$x_r = 0$ for all $r$ such that $\Psi(r) < \Psi(q)$.
Hence $\beta(q) = 0$ {\em and} $\beta(q) = x_q \Lambda_q^-$, which is impossible.
\end{proof}

This corollary leads to the following properties of canonical classes.

\begin{lemma}\label{le:canprop}
Let a torus $T$ act on a compact oriented manifold $M$ with isolated fixed
points,  and let $\Psi \colon  M \to \R$ be a formal moment map.
For all $p \in M^T$,
the canonical class $\alpha_p$ is uniquely determined by the
orientation on $\nu^-(p)$.
\end{lemma}

\begin{proof}
Fix an orientation on $\nu^-(p)$
and let $\alpha_p$  and $\alpha'_p$ be canonical classes for $p \in M^T$.
Consider the class $\beta = \alpha_p - \alpha'_p
\in H^{2 \lambda(p)}_T(M;\Q)$.
Assume $\beta \neq 0$.
Then, since the restriction map from
$H_T^*(M;\Z)$  to $H_T^*(M^T;\Z)$ is
injective, there exists $q \in M^T$ such that $\beta(q) \neq 0$
but  $\beta(r) = 0$
for all $r \in M^T$ satisfying $\Psi(r) < \Psi(q)$.
  By the definition of canonical class,
$\beta(s) = 0$ for all $s \in M^T$ such that $\lambda(s) \leq \lambda(p)$;
therefore $\lambda(q) > \lambda(p)$.
But this contradicts Corollary~\ref{co:towardindexstatement}.
\end{proof}
\begin{lemma}\label{le:2prime}
Let a torus $T$ act on a compact oriented manifold $M$ with isolated fixed
points,   and
let $\Psi \colon  M \to \R$ be a formal moment map. If
$\alpha_p\in H_T^{2\lambda(p)}(M;\Q)$ is a canonical class for a fixed point
$p$, then $\alpha_p$ also satisfies the property:
\begin{itemize}
\item[$(2^\prime)$]
$\alpha_p(q) = 0$ for all $q \in M^T \smallsetminus \{p\}$ such that $\Psi(q) \leq \Psi(p)$.
\end{itemize}
\end{lemma}

\begin{proof}
There exists a point $q \in M^T$ so that
$\alpha_p(q)\neq 0$ but $\alpha_p(r) = 0$ for all $r \in M^T$ so
that $\Psi(r) < \Psi(q)$.
By the definition of canonical classes, the fact
that $\alpha_p(q) \neq 0$ implies that either $q = p$ or
$\lambda(q)>\lambda(p)$.
In the latter case,
Corollary~\ref{co:towardindexstatement} implies that $\alpha_p(q) = 0$.
Thus $q = p$.
\end{proof}

This has the following important consequence.

\begin{corollary}\label{co:increasing}
Let a torus $T$ act on a compact oriented manifold $M$ with isolated fixed
points,  and let $\Psi \colon  M \to \R$ be a formal moment map.
Define an oriented graph with vertex set $V = M^T$ and edge
set
$$
\onestep = \left\{ \left. (r,r') \in M^T \times M^T \right| \lambda(r') - \lambda(r) = 1 \mbox{ and } \alpha_r(r') \neq 0
\right\}.
$$
If there exists a path in $(V, E)$ from a point $p$ to a different
point $q$, then
$\Psi(p) < \Psi(q)$.
\end{corollary}

Finally,  canonical classes do exist
in a number of important special cases. 
For example, the following lemma is a special case
of Lemma 1.13 in \cite{MT:fundamentalgroups}.

\begin{lemma} [McDuff-Tolman]
Let the circle $S^1$ act on a compact oriented manifold $M$ with isolated fixed points,
 and let $\Psi \colon  M \to \R$ be a formal moment map.
Then there exists a canonical class $\alpha_p \in H^*_{S^1}(M;\Q)$ for
all $p \in M^T$.
\end{lemma}

Note that if $\Psi$ is a formal moment map then
$-\Psi$ is also a formal moment
map. If the index of $\Psi$ at $q \in M^T $ is
$2 \lambda(q)$, the index of $-\Psi$ at $q$ is
$\dim(M) - 2\lambda(q)$.
Finally, the tangent bundle at $q$ is an oriented real
representation of $T$; let $\Lambda_q$ denote the product of
the weights of this representation.
As one might expect, canonical classes for $\Psi$ exist  for all points $p\in M^T$ exactly when
they exist for $-\Psi$ for all points.

\begin{lemma}\label{le:beta} Let a torus $T$ act on a compact oriented manifold $M$ with isolated fixed points and let
$\Psi \colon M \to \R$ be a formal
moment map. If  there exists a canonical class $\alpha_p\in H_T^{2\lambda(p)}(M;\Q)$ with respect to $\Psi$ for each
$p\in M^T$,
then there exists a canonical class
$\beta_q \in H^{\dim(M) - 2\lambda(q)}_T(M;\Q)$ with respect to $-\Psi$ for
each $q \in M^T.$
Moreover,   $\beta_q$ is integral for all $q\in M^T$ if and only if  $\alpha_p$ is integral for all $p\in M^T$.
\end{lemma}

\begin{proof}
To begin, fix an orientation on $\nu^+(r)$ for all $r \in M^T$,
and consider the formal moment map $-\Psi$.
We need to find a  class $\beta_q \in H_T^{\dim(M) - 2\lambda(q)}(M;\Z)$
such that
\begin{enumerate}
\item[$(\overline{1})$]  $\beta_q(q) = \Lambda^+_q$, and
\item[$(\overline{2})$]  $\beta_q(p) = 0$ for all
 $p \in M^T \smallsetminus\{q\}$
such that $\lambda(p) \geq \lambda(q)$.
\end{enumerate}

By Lemma~\ref{le:pclass}, for all $r \in M^T$
there exists a cohomology class $\gamma_r \in H_T^{\dim(M) - 2\lambda(r)}(M;\Z)$
such that
\begin{enumerate}
\item[$(\overline{1})$]  $\gamma_r(r) = \Lambda^+_r$, and
\item[$(\overline{2}')$]  $\gamma_r(p) = 0$ for all
 $p \in M^T \smallsetminus\{r\}$ so
that $\Psi(p) \geq \Psi(r)$.
\end{enumerate}

Now let $\beta_q$ be any class which satisfies $(\overline{1})$
and $(\overline{2}')$ for $q$.
Define
$$\mathcal{L}:=\{p\in M^T|\ \lambda(p)\geq \lambda(q)\mbox{ but }
\beta_q(p) \neq 0 \}.$$
If $\L = \emptyset$, we are done. Otherwise, let
$p \in \L$ be an element which maximizes $\Psi|_\L$.
Then
$$(\alpha_p \cup \beta_q)(p) =  \Lambda_p^- \beta_q(p),$$
where $\Lambda_p^-$ is the product of the weights on $\nu^-(p)$ with
respect to the orientation 
compatible with $\alpha_p$.
Now consider any fixed point $r \neq p$.
If $\lambda(r) \leq \lambda(q)$, then
$\lambda(r) \leq \lambda(p)$ and so
$\alpha_p(r) = 0$.
Similarly, if $\Psi(r) \leq  \Psi(p)$ then $\alpha_p(r) = 0$
by Lemma~\ref{le:canprop}.
On the other hand, if $\lambda(r) \geq \lambda(q)$ and   $\Psi(r) > \Psi(p)$
then $\beta_q(r) = 0$ since $p$ maximizes $\Psi|_\L$.
Thus the product
$$(\alpha_p \cup \beta_q)(r) = 0
\quad \mbox{for all }  r \in M^T \smallsetminus \{p\}.$$
We now integrate using the Atiyah-Bott-Berline-Vergne localization formula  to obtain
$$
\int_M \alpha_p \cup \beta_q =  \frac{\Lambda_p^- \beta_q(p) }{\Lambda_p}
= \pm \frac{\beta_q(p)}{\Lambda_p^+} \in H^*(BT;\Q).
$$
This expression is integral if $\alpha_p$ and $\beta_q$
are both integral.
Hence, the class
$\beta'_q = \beta_q -\left( \frac{\beta_q(p)}{\Lambda_p^+} \right) \beta_p$
satisfies
\begin{itemize}
\item $\beta'_q(p) = 0$, 
\item $\beta'_q(r) = 0$ for all $r \in M^T \smallsetminus \{p\}$
such that $\Psi(r) \geq \Psi(p)$, and 
\item in particular, $\beta'_q$ also satisfies
$(\overline{1})$ and $(\overline{2}')$ for $q$.
\end{itemize}
The result now follows by induction.
\end{proof}

\begin{remark} \rm
Assume that for each $p \in M^T$,  there exists a canonical class $\alpha_p \in H_T^{2 \lambda(p)}(M)$
with respect to $\Psi$  which is 
compatible with a chosen orientation on $\nu^-(p)$.
Let  $\beta_q \in H_T^{\dim(M) - 2 \lambda(q)}(M)$
be 
a canonical class with respect to $-\Psi$ for each
$q \in M^T$ described in  Lemma~\ref{le:beta}.  
Clearly,
$\beta_q$ can be chosen to be compatible with the orientation of $\nu^+(q)$ induced by the orientations on $M$ and on $\nu^-(q)$.  
Then, since $\Lambda_p = \Lambda_p^+ \Lambda_p^-$,
the sets $\{\alpha_p\}$ and $\{\beta_q\}$ are dual basis 
for $H_T^*(M)$ as an $H_T^*(M)$-module under the intersection
pairing. To see this, note that
if $\lambda(q)>\lambda(p)$, then $\int_M \alpha_p\beta_q=0$ by 
degree considerations. 
If $\lambda(q) \leq  \lambda(p)$ and $p\neq q$, then 
$( \alpha_p \cup \beta_q) (r) = 0$ for all $r \in M^T$ by
the definition of canonical class, and so
$\int_M\alpha_p\beta_q=0$  by the Atiyah-Bott-Berline-Vergne
localization formula.
Finally, by a similar argument, $\int_M \alpha_p \cup \beta_p = 1.$
\end{remark}

\section{Proof of Theorem~\ref{th:pathformula} }\label{se:induction}

We are now ready to prove Theorem~\ref{th:pathformula}.
Let a torus $T$ act on a compact symplectic manifold
$(M,\omega)$ with isolated fixed points and moment map
$\Phi \colon M \to \t^*$.
Let $\Psi = \Phi^\xi$ be a generic component of the moment map.
Assume that there exists a canonical class
$\alpha_p \in H_T^{2 \lambda(p)} (M;\Q)$ for  all $p \in M^T$.
Define an oriented graph with vertex set $V = M^T$ and edge set
\begin{equation}\label{edge}
\onestep = \left\{ \left. (r,r') \in M^T \times M^T \right| \lambda(r') - \lambda(r) = 1 \mbox{ and } \alpha_r(r') \neq 0
\right\}.
\end{equation}
We need to show that for all $p$ and $q$ in $M^T$,
\begin{equation}\label{path1}
\alpha_p(q) =
 \Lambda_q^-
\sum_{{\bf r} \in \Sigma_p^q} \
\prod_{i=1}^{|{\bf r}|}
\frac {\Phi(r_{i}) - \Phi(r_{i-1})} {\Phi(q)-\Phi(r_{i-1})}
\frac {\alpha_{r_{i-1}}(r_i)}{ \Lambda^-_{r_i}},
\end{equation}
where $\Sigma_p^q$ denotes the set of paths from $p$ to $q$ in $(V,\onestep)$.
\noindent
By Lemma~\ref{le:beta},  for all $q \in M^T$, there exists a class
$\beta_q \in H_T^{\dim(M) - 2\lambda(q)}(M;\Q)$ satisfying
\begin{enumerate}
\item[$(\overline{1})$]  $\beta_q(q) = \Lambda^+_q$, and
\item[$(\overline{2})$]  $\beta_q(p) = 0$ for all
 $p \in M^T \smallsetminus\{q\}$ such that
$\lambda(p) \geq \lambda(q)$.
\end{enumerate}
In our proof of \eqref{path1}, we will also show that
\begin{equation}\label{path2}
\beta_q(p) =    \Lambda_p^+
\sum_{{\bf r}\in \Sigma_p^q}
\prod_{i=1}^ {|{\bf r}|}
\frac {\Phi(r_{i}) - \Phi(r_{i-1}) }{ \Phi(p) - \Phi(r_i) }
\frac {\alpha_{r_{i-1}}(r_i)}{ \Lambda^-_{r_i}}
\end{equation}
for all fixed points $p$ and $q$.

Consider first the case that $\lambda(q)-\lambda(p) \leq 0$.
If $p \neq q$, then
$\alpha_p(q) = 0$ by
definition.
Moreover, there are no paths
from $p$ to $q$ in $(V,\onestep)$,
and so the right hand side of \eqref{path1}  vanishes, as required.
If $p=q$, then
$\alpha_p(q) = \alpha_p(p) = \Lambda^-_p$.
In this case, the right hand side of
\eqref{path1}
is a sum over one degenerate path ${\bf r}=(p)$,
and the product is the empty product,
so the total contribution is
$\Lambda^-_p$.
Thus \eqref{path1} also holds in this case.  A nearly identical argument
proves that \eqref{path2} is satisfied.

Next, suppose that
$\lambda(q)-\lambda(p)=1$. 
If $\alpha_p(q) \neq 0$, then there is one path $(p,q)$ 
from $p$ to $q$,
and the right hand side of \eqref{path1} is 
$\alpha_p(q)$. On the other hand, if
$\alpha_p(q) = 0$, then $(p,q) \not\in E$, and so the right hand side of
\eqref{path1} vanishes.
To prove \eqref{path2}, note that by the definition of canonical class,
$$ (\alpha_p \cup \beta_q)(p) = \Lambda_p^-  \beta_q(p) \quad\mbox{ and}\quad
(\alpha_p \cup \beta_q)(q) =  \Lambda_q^+ \alpha_p(q).$$
Now consider any fixed point $r$ which is not $p$ or $q$.
If $\lambda(r) \leq \lambda(p)$ then $\alpha_p(r) = 0$,
while  if $\lambda(r) \geq \lambda(q)$,
then $\beta_q(r) = 0$.
Therefore,
$$(\alpha_p \cup \beta_q)(r) = 0 \quad \mbox{for all } r \in M^T \smallsetminus \{p,q\}.$$
Since $\deg(\alpha_p \cup \beta_q) = \dim(M)  - 2 < \dim(M)$, the integral
of $\alpha_p \cup \beta_q$ over $M$ is zero. Thus, by the
Atiyah-Bott-Berline-Vergne localization
theorem,
$$ \int_M \alpha_p\cup\beta_q =  \frac{\beta_q(p)\Lambda^-_p}{\Lambda_p}
+ \frac{\alpha_p(q) \Lambda^+_q }{\Lambda_q}=0.$$
Therefore,
$$  \beta_q(p) =
- \Lambda_p^+ \frac{\alpha_p(q) }{\Lambda^-_q}.$$

Fix $k > 1$, and assume that \eqref{path1} and \eqref{path2} hold for all fixed points
$p$ and $q$  so that $\lambda(q) - \lambda(p) < k$.
Consider fixed points $p$ and $q$
so that $\lambda(q) - \lambda(p) = k$.
We will prove that \eqref{path1} and \eqref{path2} follow for this
$p$ and $q$.

Suppose first that  $\Phi(p) = \Phi(q);$ a fortiori
$\Psi(p) = \Psi(q)$.
Then the left hand sides of \eqref{path1} and \eqref{path2}
vanish by Lemma~\ref{le:2prime}.
Since there is no path from $p$ to $q$ by Lemma~\ref{co:increasing},
the right hand sides also vanish.

So assume instead that $\Phi(p) \neq \Phi(q)$.
Let $\widetilde{\omega} = \omega +  \Phi -  \Phi(p)$ be an equivariant
extension of $\omega$.
Since $\widetilde{\omega}(r) = \Phi(r)-\Phi(p)$ for all $r \in M^T$,
$$ (\alpha_p \cup \beta_q \cup \widetilde{\omega})(q) =
\Lambda_q^+ \alpha_p(q)(\Phi(q) - \Phi(p)), \quad \mbox{and} \quad
(\alpha_p \cup \beta_q \cup \widetilde{\omega})(p) = 0 .$$
Since  $k > 1$,
$\deg(\alpha_p \cup \beta_q \cup \widetilde{\omega}) =
\dim(M) - 2k + 2 < \dim(M),$
and so the integral of $\alpha_p \cup \beta_q \cup \widetilde{\omega}$
over $M$  is zero.
Therefore, by the Atiyah-Bott-Berline-Vergne localization formula,
$$
\int_M \alpha_p \cup \beta_q \cup \widetilde{\omega}
= \frac{\Phi(q) - \Phi(p)}{\Lambda_q^-} \alpha_p(q)+
\sum_{r \neq p, q}
\frac{(\alpha_p  \cup \beta_q \cup \widetilde{\omega})(r)}
{\Lambda_r} = 0.
$$
Since  $\Phi(p) \neq \Phi(q)$,
we can solve the above equation for $\alpha_p(q)$;
\begin{equation}
\label{alpha}
\alpha_p(q) = \frac{- \Lambda_q^-}{\Phi(q) - \Phi(p)}
\sum_{r \neq p, q}
\frac{(\alpha_p \cup  \beta_q \cup \widetilde{\omega})(r)}{\Lambda_r}.
\end{equation}

Consider any fixed point $r \neq p$ or $q$.
Assume first  that $\lambda(p) <  \lambda(r) < \lambda(q)$, and
let $l  = \lambda(r) - \lambda(p)$.
By the inductive assumption
\begin{align*}
\alpha_p(r) & =
 \Lambda_r^-
 \sum_{(s_0,\ldots,s_l) \in\, \Sigma_p^r}
\prod_{i=1}^l \frac {\Phi(s_{i})- \Phi(s_{i-1})}
 {\Phi(r)-\Phi(s_{i-1})}
\frac {\alpha_{s_{i-1}}(s_i)}{ \Lambda^-_{s_i}},
\quad \mbox{and} \\
\beta_q(r) & =  \Lambda_r^+
\sum_{ (s_{l},\ldots,s_{k}) \in\, \Sigma_r^q}
\prod_{i=l+1}^{k} \frac {\Phi(s_{i}) - \Phi(s_{i-1})}{ \Phi(r) - \Phi(s_i)}
\frac {\alpha_{s_{i-1}}(s_i)}{ \Lambda^-_{s_i}}.
\end{align*}
Therefore,
if  $\Sigma_p^q(r)$ denotes the set of paths in $\onestep$ from $p$ to $q$
that pass through $r$, then
\begin{align}
(\alpha_p\cup\beta_q \cup \widetilde{\omega})(r)
& =
\Lambda_r^+ \Lambda_r^- (\Phi(r) - \Phi(p)) \sum_{{\bf s}  \in \Sigma_p^q(r)}
\frac
{\prod_{i=1}^k  [\Phi(s_i)-\Phi(s_{i-1})]\frac {\alpha_{s_{i-1}}(s_i)}{ \Lambda^-_{s_i}} }
 {\prod_{i \in \{0,\ldots,k\} \smallsetminus \{l\}}  \Phi(r)-\Phi(s_i)}
\notag \\
\label{restriction}
& = \Lambda_r \sum_{{\bf s}  \in \Sigma_p^q(r)}
\prod_{i=1}^k
[\Phi(s_i)-\Phi(s_{i-1})]
\frac {\alpha_{s_{i-1}}(s_i)}{ \Lambda^-_{s_i}}
\prod_{i\neq l} \frac{1}{\Phi(r)-\Phi(s_i)},
\end{align}
where we use the expression $\prod_{i \neq l}$ as shorthand
for $\prod_{i \in \{1,\ldots,k\} \smallsetminus \{l\}}$.
On the other hand, if $\lambda(r) \leq \lambda(p)$ or
$\lambda(q) \leq \lambda(p)$,
then $(\alpha_p \cup \beta_q \cup \widetilde{\omega})(r) = 0$.
By Lemma~\ref{co:increasing}, the right hand side of
\eqref{restriction} also vanishes.
Therefore, \eqref{restriction} holds for all
$r \in M^T \smallsetminus \{p,q\}$.

Substituting \eqref{restriction} into \eqref{alpha},
we see that
\begin{align*}
\alpha_p(q) & =
\frac{- \Lambda_q^-}{\Phi(q) - \Phi(p)}
\sum_{r \in M^T \smallsetminus \{p,q\}}
\left(
\sum_{{\bf s}  \in \Sigma_p^q(r)}
\prod_{i=1}^k
[\Phi(s_i)-\Phi(s_{i-1})]
\frac {\alpha_{s_{i-1}}(s_i)}{ \Lambda^-_{s_i}}
\prod_{i\neq l}\frac{1}{ \Phi(r)-\Phi(s_i)}
\right)
 \\
& = \frac{- \Lambda_q^-}{\Phi(q) - \Phi(p)}
 \sum_{ {\bf s} \in\, \Sigma_p^q }
 \prod_{i=1}^k
[\Phi(s_i) - \Phi(s_{i-1})]
\frac {\alpha_{s_{i-1}}(s_i)}{ \Lambda^-_{s_i}}
\left( \sum_{l=1}^{k-1} \prod_{i \neq l}  \frac{1}{ \Phi(s_l) - \Phi(s_i)}
\right) \\
& = \frac{- \Lambda_q^-}{\Phi(q) - \Phi(p)}
 \sum_{ {\bf s} \in\, \Sigma_p^q }
 \prod_{i=1}^k
[\Phi(s_i) - \Phi(s_{i-1})]
\frac {\alpha_{s_{i-1}}(s_i)}{ \Lambda^-_{s_i}}
\left( - \prod_{i = 1}^{k-1}  \frac{1}{ \Phi(s_k) - \Phi(s_i)}
\right) \\
& =  \Lambda_q^-
 \sum_{ {\bf s} \in\, \Sigma_p^q }
 \prod_{i=1}^k
[\Phi(s_i) - \Phi(s_{i-1})]
\frac {\alpha_{s_{i-1}}(s_i)}{ \Lambda^-_{s_i}}
\prod_{i = 1}^{k}  \frac{1}{ \Phi(q) - \Phi(s_{i-1})},
\\
\end{align*}
where the third equality is by   Lemma~\ref{le:CPnCancellation}.
The proof of \eqref{path2} is nearly identical.

\begin{remark} \rm
In fact, the same proof works if a torus $T$ acts on a manifold $M$ with
isolated fixed points, $\Psi \colon M \to \R$ is a formal moment map,
and $\widetilde{\eta} = \eta + \Phi$ is any closed equivariant $2$-form
(not necessarily symplectic) so that $\Phi(p) \neq \Phi(q)$ for
every pair of fixed points $p$ and $q$ so that there is
a path (of length two of more) in $(V,E)$ from $p$ to $q$.
\end{remark}

The following corollary is immediate.

\begin{corollary}
Consider the situation described in Theorem~\ref{th:pathformula}.
If $\Psi = \Phi^\xi$
achieves its minimum value at $p$,
then for any fixed point $q$,
$$
1 =
\Lambda_q^-
\sum_{{\bf r} \in \Sigma_p^q} \
\prod_{i=1}^{|{\bf r}|}
\frac {\Phi(r_{i}) - \Phi(r_{i-1})} {\Phi(q)-\Phi(r_{i-1})}
\frac {\alpha_{r_{i-1}}(r_i)}{ \Lambda^-_{r_i}},
$$
where $\Sigma_p^q$ denotes the paths in $(V,E)$ from
$p$ to $q$.
In particular, every fixed point is connected by
a path in $(V,E)$ to the minimum
(and to the maximum).
\end{corollary}

\begin{proof}
Since $\lambda(p) = 0$, $\alpha_p=1$;
hence,  $\alpha_{p}(q) = 1$ for all fixed points $q$.
\end{proof}

Our proof of Theorem~\ref{th:pathformula} relies on the following fact, which was also proved in \cite{GZ:Thomclasses} using different techniques.

\begin{lemma}\label{le:CPnCancellation}
Given  $k > 1$ distinct vectors $v_1,\ldots,v_k$ in a vector space $V$,
 $$\sum_{l=1}^k \ \prod_{i \neq l} \frac{1}{v_i -v_l} = 0.$$
\end{lemma}

\begin{proof}
The $(S^1)^k$ action on $\C^k$ induces a symplectic action of
on $\CP^{k-1}$ with fixed points
$p_1,\ldots,p_k$.
so that the  weights at $p_l$
are $\{x_i - x_l\}_{i \neq l} \in H^*((B(S^1)^k;\Z) = \Z[x_1,\ldots,x_k]$.
Since $\deg(1) = 0 < 2k - 2$, the integral of $1$ over $\CP^{k-1}$ is $0$.
Therefore, by the  Atiyah-Bott-Berline-Vergne localization theorem,
the following equation holds in the field of rational functions
$\Q(x_1,\ldots,x_k)$:
$$ \int_{\CP^{k-1}} 1 =  \sum_{l=1}^k \  \prod_{i \neq l} \frac{1}{x_i - x_l} =0.$$
Since the $v_i$'s are distinct, the claim follows easily.
\end{proof}

\section{The GKM case}\label{se:GKM}

\begin{example}\label{ex:indexincreasing} Here generic component of how Example 2.1 has an index-increasing component of moment map.
I think that we can skip this!  Sue.
\end{example}

The main goal of this section is to prove Theorem~\ref{th:GKMpathformula}.
In fact, this theorem
is an immediate consequence of Theorem~\ref{th:pathformula}
and the theorem below.

\begin{theorem}\label{th:GKMonestep}
Let $(M,\omega,\Phi)$ be  GKM space, and let
$(V,\GKMedges)$ be the associated GKM graph.  Let
$\Psi = \Phi^\xi$ be
a generic component
of the moment map; assume
that $\Psi$ is index increasing.
Then
\begin{itemize}
\item There exists a canonical class
$\alpha_p\in H_T^{2 \lambda(p)}(M;\Z)$ for all $p \in M^T$.
\item
Given $p$ and $q$ in $M^T$ such that $\lambda(q) - \lambda(p) = 1$,
\begin{equation}
\alpha_p(q) =
\begin{cases}
\Lambda_q^- \frac{\Theta(p,q)}{\eta(p,q)}  & \mbox{if } (p,q) \in \GKMedges,
 \mbox{ and} \\
0 & \mbox{if } (p,q) \not\in \GKMedges.
\end{cases}
\end{equation}
\item $\Theta(r,r')\in \Z \smallsetminus \{0\}$ for all $(r,r')\in \GKMedges$.
\end{itemize}
\end{theorem}

\begin{remark}\em
Conversely,  if there exists a canonical class $\alpha_p \in H_T^*(M;\Q)$  
for all $p \in M^T$,
then $\Psi$ is index increasing.
To see this, suppose that there
exists an edge $(p,q) \in \GKMedges$ so that $\Psi(p) < \Psi(q)$ and $\lambda(p) \geq \lambda(q)$.
On the one hand, since $\alpha_p$ is a canonical class this implies that $\alpha_p(q) = 0$.
On the other hand,  compatibility along $(p,q)$ guarantees that
$\alpha_p(q) - \alpha_p(p)$ is a multiple of $\eta(p,q)$.
Sine $\alpha_p(p) = \Lambda_p^-$ is not a multiple of $\eta(p,q)$, this implies that $\alpha_p(q) \neq 0$.
\end{remark}

\begin{remark}\label{rm:increasing}\em
There are several situations where we can immediately conclude
that $\Psi$ is index increasing.

For example, if there is an $T$-invariant Palais-Smale metric $g$ on $M$,
then $\Psi$ is index increasing.
To see this, consider an edge $(p,q) \in \GKMedges$ such
that $\Psi(p) < \Psi(q)$.
There exists a codimension one subgroup $K \subset T$
so that $p$ and $q$ are contained in the same component
$N$ of $M^K$.
Since the metric is $T$-invariant,  
$N  \smallsetminus \{p,q\}$
must be contained in both the flow up from
$p$ and the flow down from $q$. Since these flows intersect
transversally, this implies that 
the intersection has dimension $2\lambda(q)-2\lambda(p)\geq 2$, which proves $\lambda(p) < \lambda(q)$.

Similarly, if $H^{2i}(M;\Q) = \Q$ for all $i$ such
that $0 \leq 2i \leq \dim(M)$, then 
Proposition 3.4 in \cite{T} implies that every 
generic component of the moment map is index increasing.
\end{remark}

Our proof of the theorem above relies  heavily on a technical proposition,
Proposition~\ref{pr:GKM}, for which
we need a few definitions.
Let $(M,\omega,\Phi)$ be a GKM space,
and let $(V,\GKMedges)$ be the associated GKM graph.
Fix a generic
$\xi \in \t$ and consider the Morse function $\Psi = \Phi^\xi.$

A path $r$ in $(V,\GKMedges)$ is {\bf ascending} if $\Psi(r_i) \leq \Psi(r_{i+1})$ for all $i$ such that $0 \leq i \leq |r|$;
it is {\bf descending} if $\Psi(r_i) \geq \Psi(r_{i+1})$ for all such $i$.
Given $p \in M^T$, the {\bf stable set} of $p$, denoted $V_p$,   is the set of $q \in V$ such that there exists an
ascending  path  from $p$ to $q$ in $(V,\GKMedges)$;
the {\bf unstable set} of $p$, denoted $V^p$,   is the set of $q \in V$ such that there exists an
descending  path  from $p$ to $q$.
Note that $p$ itself lies in both $V_p$ and $V^p$.
Moreover, since $(r,r') \in \GKMedges$ exactly if $(r',r) \in \GKMedges$,
$p \in V_q$ exactly if $q \in V^p$.

\begin{prop}\label{pr:GKM}
Let $(M,\omega,\Phi)$ be a GKM space, and let $\Psi = \Phi^\xi$ be a generic
component of the moment map.
\begin{itemize}
\item[(a)]
For every fixed point $p$, there exists a
class $\alpha_p \in H_T^{2\lambda(p)}(M;\Z)$ such that
\begin{itemize}
\item [$(1)$]$\alpha_p(p) = \Lambda_p^-$, and
\item  [$(2^{\prime\prime})$] $\alpha_p(q)=0$ for all $q \in M^T \smallsetminus V_p$.
\end{itemize}
\item [(b)]
Given  a class $\beta \in H_T^*(M;\Z)$ and point  $q \in M^T$,
\begin{equation}\label{eq:betaproduct}
\beta(q) =  x\spaceafterx \prod_{\substack{
(r,q) \in \GKMedges \\ \beta|_{V^r} = 0 }}
\spacebeforeeta \eta(r,q),
\quad \mbox{where } x \in H^*(BT;\Z).
\end{equation}
\end{itemize}
\end{prop}

Proposition \ref{pr:GKM}(a) is proved for  {\em rational} classes in
the more general setting of GKM graphs in \cite{GZ:Thomclasses}.

\begin{remark}\rm
Part (a) of this proposition is exactly what geometric intuition leads you to expect.
To see this, fix a generic $T$-invariant metric on $M$.
As the name suggests, the stable set $V_p$ should  be the set of vertices which
are in the closure of the stable manifold of $p$.
Moreover,  one should be able to adapt the Morse theoretic proof of Lemma~\ref{le:pclass}
to directly prove that the class $\alpha_p$ which you construct is supported  on the set of  vertices
in the closure of the stable manifold of $p$.
\end{remark}

\begin{remark}
\rm
Part (b) is slightly more subtle.
Given a class $\beta \in H_T^*(M;\Q)$ and an edge $(r,r') \in \GKMedges$,
$\beta(r) - \beta(r')$ is  a rational multiple\footnote{
In fact, a  rational class $\beta\in H_T^*(M^T;\Q)$
is in the image of the restriction map $\iota^*: H_T^*(M;\Q)\rightarrow H_T^*(M^T;\Q)$
exactly if $\beta(q)-\beta(p)$  is a multiple of $\eta(p,q)$ for every
edge $(p,q)\in \GKMedges$ (see \cite{GKM}).}
of $\eta(r,r')$.
Since the weights at each fixed point are pairwise
linearly independent, this immediately implies that for any fixed point $q$,
$$
\beta(q) = x\spaceafterx\prod_{\substack{(r,q)\in \GKMedges \\ \beta(r)=0}} \spacebeforeeta\eta(r,q),
\quad \mbox{where } x \in H^*(BT;\Q).
$$
However, the same statement is not true integrally.
Although $\beta(q)$ must be an integral multiple of $\eta(r,q)$ for each $(r.q) \in \GKMedges$
such that $\eta(r) = 0$, it might not be an integral multiple of the product of
these weights because the weights  might not be pairwise relatively prime; see
Example~\ref{ex:notprime}.
Notice that expression \eqref{eq:betaproduct} has fewer terms in the product.
\end{remark}

\begin{example}\label{ex:notprime}\rm
Let $S^1 \times S^1$ act on $\CP^1 \times \CP^1$ by
$$(t_1,t_2) \cdot \left( [z_0: z_1], [w_0:w_1] \right) =
\left( [t_1^2\ z_0: z_1], [t_2^2\ w_0:w_1] \right);$$
let $\Phi \colon M \to \t^*$ be the moment map.
Let $\xi = (1,1) \in \R^2$ and let $\Psi = \Phi^\xi$.
The associated GKM graph has four vertices,
$$ SS = ([1:0],[1:0]), \ SN = ([1:0],[0:1]), \ NS = ([0:1],[1:0]), \ \mbox{and }
NN = ([0:1],[0:1]).$$
There are four ascending edges,
$(SS,NS)$,  $(SN, NN)$,
$(SS, SN)$,  and $(NS,NN)$ (and hence also four descending edges).
The first two have weight $2 x_1$, the latter two have weight $2 x_2$.
There exists a class $\beta \in H_T^4(M;\Z)$ so that
$$\beta(SN) = \beta(NS) = 0, \quad \mbox{and} \quad
\beta(SS) = \beta(NN) = 2 x_1 x_2.$$
Although $\beta(SN)  = 0$,
$\beta|_{V^{SN}} \neq 0 $ because
$V^{SN} = \{SN, SS\}$.
Similarly, $\beta|_{V^{SN}} \neq 0.$
Therefore, even though $\beta_{NN} = 2 x_1 x_2$ is not a multiple of $\Lambda_{NN}^- = \eta(SN,NN) \cdot \eta(NS,NN)
= 4x_1 x_2$, this example
 does satisfy~\eqref{eq:betaproduct} in Proposition~\ref{pr:GKM}.
\end{example}

We begin by showing that Theorem~\ref{th:GKMonestep} follows from Proposition~\ref{pr:GKM}.
\begin{proof}[Proof of Theorem~\ref{th:GKMonestep}]
Fix $p \in M^T$. By Proposition~\ref{pr:GKM}, there exists a class
$\alpha_p \in H_T^{2 \lambda(p)}(M;\Z)$ which satisfies properties
$(1)$ and $(2^{\prime\prime})$.  Since $\Psi$ is index increasing, $\lambda(q) > \lambda(p)$
for all $q  \in V_p \smallsetminus \{p\}$.   Hence,
$\alpha_p$ is a canonical class.

Consider $q \in M^T$ such that $\lambda(q) = \lambda(p) + 1$.
If $(p,q) \not\in \GKMedges$, then $q \not\in V_p$, and
so $\alpha_p(q) = 0$.  Now assume
that $(p,q) \in \GKMedges$.
There are  $\lambda(p) = \lambda(q) - 1$
{\em other} edges $e_1,\ldots, e_{\lambda(p)}\in\GKMedges$
of the form  $e_i = (p_i,q)$
with $\Psi(p_i)< \Psi(q)$.
Since $\Psi$ is index increasing,
$\lambda(p_i)\leq\lambda(p)$ and so
$V_p \cap V^{p_i} = \emptyset$
for all $i$.
Thus  $\alpha_p|_{V^{p_i}} = 0$
for all $i$.
By Proposition~\ref{pr:GKM},  $\alpha_p(q)$ is an integral multiple
of
$\prod_{i} \eta(p_i,q)= \frac{\Lambda_q^-}{\eta(p,q)}$; by degree considerations this implies
that
$$
\alpha_p(q) =\theta \cdot\frac{\Lambda_q^-}{\eta(p,q)},\quad \mbox{where }\theta \in \Z.
$$
On the other hand, compatibility
along $(p,q)$ guarantees that
$\alpha_p(q) - \alpha_p(p) = \theta \cdot \Lambda_q^-/\eta(p,q) -  \Lambda_p^-$
is a multiple of $\eta(p,q)$, and hence
$\theta$ satisfies
\begin{equation}\label{eq:integralthetadeterminedby}
 \frac{\Lambda_q^-}{\eta(p,q)}\theta = \Lambda_p^- \mod \eta(p,q).
\end{equation}
Since
$\Lambda_p^-$ is not a multiple of $\eta(p,q)$, the integer $\theta$ is
nonzero.
Finally,  a straightforward calculation shows that
\begin{equation}
\theta = \frac{\rho_{\eta(r,r')} \left(\Lambda^-_r \right)}
{\rho_{\eta(r,r')} \left(\frac{\Lambda^-_{r'}}{\eta(r,r')}\right)}.
\end{equation}
Hence, $\theta = \Theta(r,r')$.
\end{proof}

We now turn to the proof of Proposition~\ref{pr:GKM}.
We begin by establishing some terminology, after which we prove two key lemmas. Finally we use these lemmas and an inductive argument on the dimension of $M$ to complete the proof.

\begin{definition}\label{de:propQ}
Let $(M,\omega,\Phi)$ be a GKM space, and let $\Psi = \Phi^\xi$ be a generic
component of the moment map.
We say that $\beta\in H_T^*(M;\Z)$ is {\bf robustly zero}
 at $r\in M^T$ if $\beta|_{V^r} =0$. We say that $\beta$ is {\bf robustly integral}
at $q\in M^T$ if
$$
\beta(q) = x\spaceafterx\prod_{\substack{(r,q)\in \GKMedges \\ \beta|_{V^r}=0}}\spacebeforeeta\eta(r,q),\quad\mbox{where } x\in H^*(BT;\Z).
$$
\end{definition}

Thus Proposition~\ref{pr:GKM}(b) is the statement that every class
$\beta \in H_T^*(M;\Z)$
is robustly integral at every point $q\in M^T$.
We
first show that Proposition~\ref{pr:GKM}(b) implies \ref{pr:GKM}(a).

\begin{lemma}\label{le:GKM1b}
Let $(M,\omega,\Phi)$ be a GKM space, and let $\Psi = \Phi^\xi$ be a generic
component of the moment map.
Assume that every class $\beta \in H^*_T(M;\Z)$ is robustly integral at every $q \in M^T$.
Then, for every fixed point $p$, there exists a  class $\alpha_p \in H_T^{2\lambda(p)}(M;\Z)$ such that
\begin{enumerate}
\item[$(1)$] $\alpha_p(p) = \Lambda_p^-$, and
\item  [$(2^{\prime\prime})$]$\alpha_p(q)=0$ for all $q \in M^T \smallsetminus V_p$.
\end{enumerate}
\end{lemma}
\begin{proof}
By Lemma~\ref{le:pclass}, for all $q \in M^T$ there   exists a class
$\gamma_q \in H_T^{2 \lambda(q)}(M;\Z)$
\begin{enumerate}
\item [$(1)$] $\gamma_q(q) = \Lambda_q^-$.
\item [$(2^\prime)$] $\gamma_q(r) = 0$ for all $r \in M^T \smallsetminus \{q\}$
such that $\Psi(r) \leq \Psi(q)$.
\end{enumerate}
Let $\alpha_p \in H_T^{2 \lambda(p)}(M;\Z)$
be any class which satisfies $(1)$ and $(2^\prime)$ for $p$.
Define $$ \mathcal{L} = \{q\in M^T\smallsetminus V_p\mid\alpha_p(q)\neq 0\}.$$
If $\L = \emptyset$, we are done.  Otherwise, let
$q \in \L$ be an element which minimizes $\Psi|_{\L}$.
Consider any $r \in V^q \smallsetminus \{q\}$.
Since
$\Psi(r) < \Psi(q)$ and
$q$ is minimal, $r\not\in\mathcal{L}$.
Moreover, since $q \in V_r$  but $q \not\in V_p$, we
conclude that
$r\in M^T\smallsetminus V_p$; hence $\alpha_p(r)=0$.
Therefore, $\alpha_p|_{V^r}=0$ for each edge $(r,q)\in \GKMedges$ such that $\Psi(r)<\Psi(q)$.
Since $\alpha_p$ is robustly integral at $q$
this implies that
there exists $x \in H^*(BT;\Z)$ so that
$\alpha_p(q) = x \Lambda_q^- = x \gamma_q(q)$.
The difference  $\alpha'_p = \alpha_p - x \gamma_q$ satisfies
\begin{itemize}
\item $\alpha'_p(q) = 0$, 
\item $\alpha'_p(r) = \alpha_p(r)$ for all $r \in M^T \smallsetminus \{q\}$
such that $\Psi(r) \leq \Psi(q)$, and 
\item in particular, $\alpha'_p$ also  satisfies $(1)$ and $(2^{\prime})$ for $p$.
\end{itemize}
The result now follows by induction.
\end{proof}

Now we show that Proposition~\ref{pr:GKM}(a) plus the assumption that all the $\alpha_p$ are robustly integral
at all fixed points implies Proposition~\ref{pr:GKM}(b).

\begin{lemma}\label{le:GKM2b}
Let $(M,\omega,\Phi)$ be a GKM space, and let $\Psi = \Phi^\xi$ be a generic
component of the moment map.
Assume that, for every fixed point $p$, there exists a
class $\alpha_p \in H_T^{2\lambda(p)}(M;\Z)$
such that
\begin{enumerate}
\item[$(1)$] $\alpha_p(p) = \Lambda_p^-$, and
\item  [$(2^{\prime\prime})$]$\alpha_p(q)=0$ for all $q \in M^T \smallsetminus V_p$.
\item [$(3)$] $\alpha_p$ is robustly integral at every point $q\in M^T$.
\end{enumerate}
Then every class $\beta \in H_T^*(M;\Z)$ is robustly integral at every point $q\in M^T$.
\end{lemma}

\begin{proof}
Fix a class $\beta \in H^*_T(M;\Z)$.
By Lemma~\ref{le:pclass}, we can write
$$\beta = \sum_{p \in M^T} x_p \alpha_p, \quad \mbox{where } x_p \in H^*(BT;\Z)
\ \forall \  p \in M^T.$$

Fix $r \in M^T$ such that $\beta|_{V^r} = 0$.
We claim that $x_p = 0$ for all $p \in V^r$.
If not, then  there exists $p \in V^r$ so that $x_p \neq 0$ but $x_q = 0$ for all
$q \in V^r$ such that $\Psi(q) < \Psi(p)$.  But then $\beta(q) = x_q
\Lambda_q^-\neq 0$, which
contradicts the assumption that $\beta|_{V^r} = 0$.
Therefore, if $x_p \neq 0$, then $p \not\in V^r$, and
so $V_p \cap V^r = \emptyset$; hence $\alpha_p|_{V^r} = 0$.
Since each $\alpha_p$ is
robustly integral at every point $q \in M^T$, this completes
the proof.
\end{proof}

\begin{proof}[Proof of Proposition~\ref{pr:GKM}.]

Assume that Proposition~\ref{pr:GKM}
is true for all manifolds of dimension less than $\dim M$; we
will prove that it is true for $M$.
The result will then follow by induction.

We first consider the case when $T$ acts on $M$ effectively.
Choose any  $\beta\in H_T^*(M;\Z)$ and fixed point $q\in M^T$.
By Lemma~\ref{le:GKM1b}, it is enough to prove that
$\beta$  is robustly integral at $q$.

The proposition is obvious when $\dim M=2$.
Hence, $\beta(q)$ is an integer multiple
of $\eta(r,q)$  for every $(r,q) \in \GKMedges$ 
with
$\beta|_{V^r}=0$.
Since these weights are
pairwise linearly independent,
this implies that $\beta(q)$ is a
{\em rational} multiple of
\begin{equation}\label{eq:productstronglyzero}
\prod_{\substack{(r,q)\in \GKMedges\\ \beta|_{V^r}=0}}\spacebeforeeta \eta(r,q).
\end{equation}

Now consider any prime $k$ and natural number $l$ such that  $k^l$ divides the
product~\eqref{eq:productstronglyzero}.
Let $\Gamma = \{ t \in T \mid t^k = 1 \},$
and let  $N \subset M$ be the component of $M^\Gamma$ containing $q$.
Then $(N,\omega|_N,\Phi|_N)$ is also a GKM space, and $\Psi|_N =
\Phi^\xi|_N$
is a generic component of the moment map.  Let $(V_N,(\GKMedges)_N)$
be the associated
GKM graph. Given any  edge $(r,q) \in \GKMedges$,
the weight $\eta(r,q)$ is a multiple of $k$ exactly if
$(r,q) \in (\GKMedges)_N$.
Therefore, since 
$k^l$ divides the product~\eqref{eq:productstronglyzero},
it divides the (smaller) product
\begin{equation}\label{eq:restrictedproduct}
\prod_{\substack{(r,q) \in (\GKMedges)_N\\ \beta|_{V^r}=0}}\spacebeforeeta \eta(r,q).
\end{equation}
Since
$(\GKMedges)_N \subset \GKMedges$,
if $\beta\in H_T^*(M;\Z)$ is robustly zero at $r\in V_N$,
then $\iota_N^*\beta \in H_T^*(N;\Z)$ is robustly zero at $r$,
where  $\iota_N^*:H_T^*(M)\rightarrow H_T^*(N)$ is the restriction map.
Moreover, since
the $T$ action on $M$ is effective, $\dim N < \dim M$.
By the inductive hypothesis,
this implies that $\iota_N^*\beta(q)=\beta(q)$ is an integral
multiple of the product \eqref{eq:restrictedproduct}.
Therefore, since
$k^l$ divides this product,
$k^l$ also divides $\beta(q)$.
This proves the proposition when $T$ acts effectively.

We now consider the general case.
By Lemma~\ref{le:GKM2b}, it is enough to show that for each $p$ there exists
a class $\alpha_p$
which satisfies $(1)$ and $(2^{\prime\prime})$
and is robustly integral at every point $q\in M^T$.

Let
$$ \Stab(M)  = \{t \in T \mid t \cdot m = m\ \mbox{for all }m\in M \}, $$
and $T' = T/\Stab(M)$.
 Let $\Pi \colon T \to T'$
be the natural projection, and
let $\pi \colon \t \to \t'$ be the induced map on the lie algebras.
Notice that $\pi$ takes the lattice $\ell \subset \t$ to
a sublattice of $\ell' \subset \t'$.

Since $T'$ acts naturally on $M$, the
  map $\Pi \colon T \to T'$
induces  maps in equivariant cohomology%
\footnote{To see this, let $ET$ and $ET'$ be contractible spaces on
which $T$ and $T'$, respectively,  act freely. Then $T$ acts on
$ET \times ET'$ by $t \cdot (e, e') = (t \cdot e, \pi(t)\cdot e')$.
The projection from $M \times (ET \times ET')$ to $M \times ET'$
descends to a map from $M \times_T (ET \times ET')$
to $M \times_{T'} ET'$. This induces a map from
$H^*(M \times_{T'} ET') = H_{T'}^*(M)$
to $H^*(M \times_T (ET \times ET')) = H_T^*(M)$.}
$$
\Pi_M^*: H_{T'}^*(M;\Z) \longrightarrow H_T^*(M;\Z) \quad \mbox{and}  \quad
\Pi_p^*: H_{T'}^*(\{p\};\Z) \longrightarrow H_T^*(\{p\};\Z)  \ \forall \ p \in M^T.
$$
Since $M^T = M^{T'}$,  these  fit together into the following commutative diagram.
$$
\begin{CD}
 H_{T'}^*(M;\Z) @>>> \bigoplus_{p} H_{T'}^*(\{p\};\Z) \\
  @VV{\Pi_M^*}V                      @VV \bigoplus_p \Pi_p^*V\\
  H_{T}^*(M;\Z)   @>>>  \bigoplus_{p} H_T^*(\{p\};\Z).
\end{CD}
$$
Moreover, for each $p \in M^T$, 
in degree two  the map  $\Pi_p^*$ 
is identified
with the dual map $\pi^*:(\t')^*\rightarrow \t^*$ under the natural
identification of $H_T^2(\{p\};\Z)$ and $H_{T'}^2(\{p\};\Z)$ with the
lattices $\ell^*$ and $(\ell')^*$,
respectively.

Clearly, there exists a moment map $\Phi' \colon M \to (\t')^*$
so that $(M,\omega,\Phi')$ is a GKM space.
Since $\Phi' \circ \pi^* = \Phi$,
$\Psi = \Phi^\xi = (\Phi')^{\pi(\xi)}$
 is  a generic component of the moment map.
Note that the vertices and edges of the graphs for the
$T$ action and the $T'$ action are naturally identical,
and that $\pi^*$
takes the $T'$-weight of each edge to the $T$-weight of
the same edge.
In particular, $\pi^*\left((\Lambda')^-_p \right) = \Lambda_p^-$,
where  $(\Lambda')^-_p$ denotes the $T'$-equivariant Euler class
of the negative normal bundle of $\Psi$ at $p$.

Since $T'$ acts effectively on $M$, 
the first part of this proof
implies there exists a class $\alpha'_p \in H_{T'}^{2 \lambda(p)}(M;\Z)$ satisfying
\begin{enumerate}
\item[$(1)$] $\alpha'_p(p) = (\Lambda')^-_p$.
\item [$(2^{\prime\prime})$] $\alpha'_p(q) = 0$ for all $q \in M^T\smallsetminus V_p$
\item [$(3)$] $\alpha'_p$ is robustly integral 
at every point $q \in M^T$.
\end{enumerate}
Let $\alpha_p:= \Pi^*_M(\alpha'_p)\in
 H_T^*(M;\Z)$. Then $\alpha_p$ has the desired properties.
\end{proof}

\section{Examples}\label{se:examples}

We conclude our paper with two explicit examples.
\begin{example} \rm
$Fl(\C^n)$.
The standard action of $T = (S^1)^n$ induces an action on
the complete flag manifold $Fl(\C^n)$
whose fixed points are
given by flags in the coordinate lines $\C^n$. These are
indexed by permutations $\sigma\in S_n$ on $n$ letters;
the fixed point corresponding to $\sigma$ is given by
$$
\langle 0\rangle\subset \langle f_{\sigma(1)}\rangle \subset \langle f_{\sigma(1)},f_{\sigma(2)}\rangle \subset \cdots \subset \langle f_{\sigma(1)},f_{\sigma(2)},\cdots f_{\sigma(n)}\rangle = \C^n,
$$
where the brackets indicate the span of the vectors
and  $f_1,\dots, f_n$ is the standard  basis of $\C^n$.
By abuse of notation, we will denote both the permutation
and the corresponding fixed point by $\sigma$.
Let $\ell(\sigma)$ denote the length of $\sigma$,
and choose a generic $\xi$ such that $\ell(\sigma)<\ell(\sigma')$
implies $\Psi(\sigma)<\Psi(\sigma')$ for all $\sigma, \sigma'\in S_n$.
Note that $(\sigma,\sigma')\in \onestep$ if there exists a
transposition $t$ such that $t\sigma=\sigma'$ and
$\ell(\sigma)=\ell(\sigma')-1$ (in contrast, $\GKMedges$
consists of all $(\sigma, \sigma')$ such that $t\sigma=\sigma')$
for some $t$). If $t=t_{ij}$ is the transposition switching $i$
and $j$ with $i<j$, then $\eta(\sigma,\sigma') = x_i-x_j$ and
by our convention, this is considered a {\em positive} weight.

We begin by showing that for all edges $(\sigma,\sigma')\in \onestep$,
$\Theta(\sigma,\sigma')=1$.  Recall that the weights that occur in
the representation of $T$ on the negative normal bundle at a point
$\sigma$ are positive weights. The length
$\ell(\sigma)$ is also the number of positive weights at $\sigma$.
The weight $\eta\neq \eta(\sigma,\sigma')$ is a positive weight at
$\sigma'$ if and only if $t_{ij}\eta$ is a positive weight $\sigma$.
Thus there is a bijection of positive weights at $\sigma$ and
positive weights at $\sigma'$ excepting $\eta(\sigma',\sigma)$.
Moreover, for each weight $\eta$ at $\sigma$, the weight
$t_{ij}\eta$ at $\sigma'$ has the property that
$\eta=t_{ij}\eta \mod \eta(\sigma,\sigma')$. It follows
that $\Theta_{\eta(\sigma,\sigma')} =1$.

Now consider any fixed point $\sigma$, and let
$\alpha_\sigma$ be the associated canonical class.
Theorem~\ref{th:GKMpathformula} says for any $\mu\in M^T$,
\begin{equation}\label{eq:flagcase}
\alpha_\sigma(\mu) = \Lambda^-_\mu \sum_{{\bf r}\in \Sigma_{\sigma}^\mu} \prod_{i=1}^{|{\bf r}|}\frac{\Phi(r_i)-\Phi(r_{i-1})}{\Phi(\mu)-\Phi(r_{i-1})}\cdot \frac{1}{\eta(r_{i-1},r_{i})}.
\end{equation}
Note that $\eta(r_{i-1},r_{i})$ is positive.
Similarly,  $\Lambda_\mu^-$ is a product of positive weights.
Finally,  $\Phi(r_i)-\Phi(r_{i-1})$ and $\Phi(\mu)-\Phi(r_{i-1})$
are positive for each $i$. More precisely,
$\Phi(r_i)-\Phi(r_{i-1}) =k_i\eta(r_{i-1},r_{i})$
where $k_i\in \Q^+$ and $\Phi(\mu)-\Phi(r_{i-1}) = \sum_{j=1}^m k_j \eta(s_{j-1},s_{j})$,
where $j$ indexes the vertices in a path
$(s_0,\dots, s_m)$ from $s_0=r_{i-1}$ to $s_m = \mu$.
Thus every term in the expression \eqref{eq:flagcase} is positive.
\end{example}

\begin{example}\label{ex:nonkahler} \rm
The integer $\Theta(p,q)$ is not always positive; when it is negative, the canonical class $\alpha_p$ restricts to a negative value at $q$. In \cite{To:nonkahler}, the second author demonstrated the existence of a GKM space that does not have a $T$-invariant K\"ahler metric. The corresponding GKM graph $\GKMedges$ can be expressed as the image of the singular set under the moment map $\Phi$, pictured in Figure~\ref{fi:nonkahler},
where we have represented each pair of edges $(p,q)$ and $(q,p)$ by one drawn edge.
\psfrag{rho}{{\ensuremath $\rho_{\eta(p,q)}$}}
\psfrag{xi}{$\xi$}
\psfrag{p}{$p$}
\psfrag{q}{$q$}
\begin{figure}[h]
\includegraphics[width=2.5in]{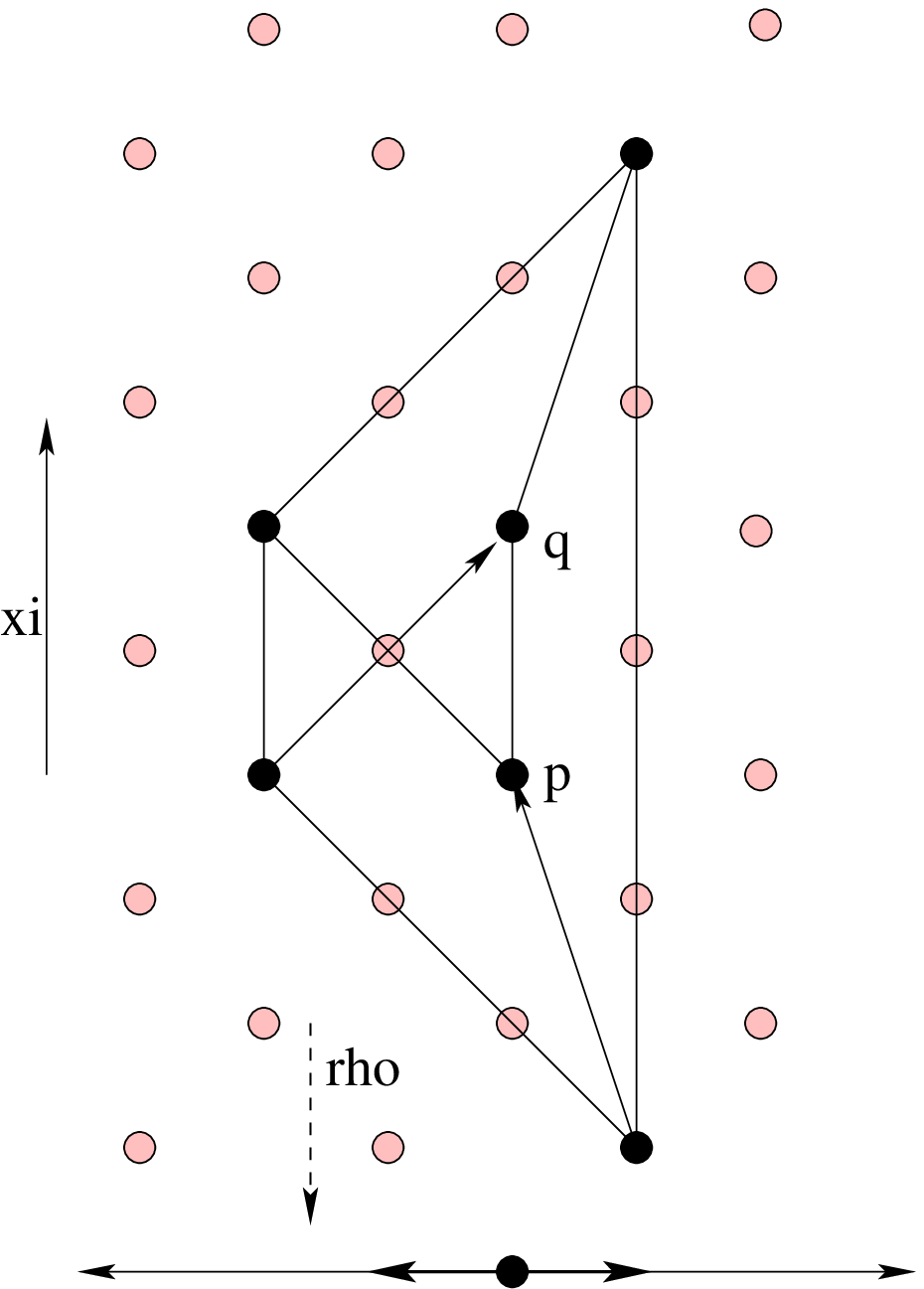}
\caption{A non-K\"ahler GKM space. The restriction of the canonical class $\alpha_p$ to $q$ is negative.}
\label{fi:nonkahler}
\end{figure}
Let $\xi\in \t$ be as indicated, and note that $\Psi=\Phi^\xi$ is an index increasing component of the moment map.
Consider the edge corresponding to $(p,q)$. Indicated on the figure are the positive weights (excluding $\eta(p,q)$) of the $T$ action on the negative normal bundles at $p$ and $q$, according to the choice of $\xi$ indicated. Under the map $\rho_{\eta(p,q)}$, these vectors project to vectors opposite in sign (and of equal magnitude).
Thus $\Theta(p,q)<0$. It immediately follows from Theorem~\ref{th:GKMpathformula} that $\alpha_p(q)<0$.

\end{example}

\end{document}